\documentclass[11pt]{article}
\newcommand{\ttau}{\tilde{\tau}}

\newcommand{\C}{\mathcal{C}}
\newcommand{\sone}{\mathbb{S}^{1}}

\newcommand{\be}{\begin{eqnarray}}
\newcommand{\ee}{\end{eqnarray}}

\newcommand{\om}{\Omega}

\newcommand{\rank}{{\rm rank}\,}

\newcommand{\1}{{\bf 1}}
\newcommand{\tr}{{\rm tr}\,}

\newcommand{\scc}{\mathcal{C}}

\newcommand{\eps}{\epsilon}

\usepackage{amsfonts}
\usepackage{amsmath}
\usepackage{amsthm}
\usepackage{epsfig}
\usepackage{amssymb}
\date{\today}
\newtheorem{theorem}{Theorem}[section]
\newtheorem{prop}{Proposition}[section]
\newtheorem{lemma}{Lemma}[section]
\newtheorem{remark}{Remark}[section]
\numberwithin{equation}{section}
\title{On one-homogeneous solutions to elliptic systems with spatial variable dependence in two dimensions}
\author{J. J.  Bevan \footnote{Address for correspondence: Department of Mathematics, University of Surrey, Guildford, Surrey, GU2 7XH, UK}}

\date{\today}
\begin{document}
\maketitle
\begin{abstract} \noindent We extend the result \cite{Ph02} of Phillips by showing that one-homogeneous solutions of certain elliptic systems in divergence form either do not exist or must be affine.  The result is novel in two ways.  Firstly, the system is allowed to depend (in a sufficiently smooth way) on the spatial variable $x$.  Secondly, Phillips's original result is shown to apply to $W^{1,2}$ one-homogeneous solutions, from which his treatment of Lipschitz solutions follows as a special case.    A singular one-homogeneous solution to an elliptic system violating the hypotheses of the main theorem is constructed using a variational method.   
\end{abstract}

\section{Introduction}

One of the main results of this paper is a regularity theorem which extends an earlier result of Phillips \cite{Ph02}.  It turns out that the extension is reasonably straightforward.  It can be used to rule out the possibility that a non-trivial (i.e., non-affine) one-homogeneous function can be a stationary point of a functional such as 
\begin{equation}\label{i}
 I(u)=\int_{\om} f(x, \nabla u(x)) \,dx,
\end{equation}
where $f$ is strongly rank-one convex in the gradient argument, sufficiently regular in the spatial variable $x$ in a two-dimensional domain $\om$ containing zero as an interior point.  Here, $u: \om \to \mathbb{R}^{m}$, where $m \geq 1$.     The second and more substantial part of the paper is devoted to finding circumstances under which the extended version of Phillips's result fails.    This involves proving the existence of a non-trivial one-homogeneous solution to an elliptic system which violates at least one of the hypotheses of the theorem.    The resulting function is $W^{1,2}$ but not Lipschitz, and it has an interesting topological effect on its domain of definition.  It is not a solution of the Euler-Lagrange equation associated with the functional it minimizes, but it does solve the so-called Equilibrium equations.  
See sections \ref{francois} - \ref{vivaldi} for more details.

Several examples of non-smooth minimizers of functionals such as \eqref{i}, though with an $x-$independent integrand,  are based on positively one-homogeneous functions: see  \cite{Ne75},  
 \cite{SY00}, \cite{SY02}. This is the main reason for our interest in one-homogeneous solutions of elliptic systems. For an overview of the example \cite{Ne75} and other singular solutions to elliptic systems see \cite{Gi83}.   It is perhaps worth pointing out that De Giorgi's example \cite{DG68} of a singular minimizer is set in three space dimensions and is based on a functional of the form \eqref{i}, although the minimizer he constructs is not one-homogeneous.   
 
 All examples cited here are set in dimensions strictly larger than two.  In two and higher dimensions \cite{MS03} and \cite{Sz04} have shown that stationary points can in general be nowhere $C^{1}$.  However, the question of the regularity of minimizers in two dimensions is still open.  This is the motivation for our study of the functional \eqref{i} above.  In order to rule out certain classes of singular solutions, such as the one-homogeneous solutions considered here, it is sufficient to show that they cannot be stationary points of $\eqref{i}$.    Together with \cite{Ph02}, this argument can be applied to the functionals appearing in \cite{Be05}, with the result that none is strongly elliptic.  

We remark that the theorem in this note and that of \cite{Ph02} are also of interest because they yield smoothness in a case not covered by the regularity theory of elliptic manifolds \cite{Sv93}. This is possibly why  one-homogeneous functions often feature in counterexamples: their tangent spaces are suitably `degenerate'.  
%always contain rank-one matrices.

\subsection{Notation and definitions}

Recall that a function $u: \mathbb{R}^{n} \to \mathbb{R}^{m}$ is positively one-homogeneous (henceforth one-homogeneous) if, for each $x$ in $\mathbb{R}^{n}$,
\begin{equation}\label{onehomdef}u(\lambda x)=\lambda u(x) \ \textrm{for all} \ \lambda \geq 0.
\end{equation}
It follows that any one-homogeneous function: $\mathbb{R}^{n} \to \mathbb{R}^{m}$ can be represented as
\begin{equation}\label{ug}u(x)= Rg\left(\frac{x}{|x|}\right),\end{equation}
where $g: \mathbb{S}^{n-1} \to \mathbb{R}^{m}$, $R = |x|$ and $\mathbb{S}^{n-1}$ is the boundary of the unit ball in $\mathbb{R}^{n}$.   
Conversely, any choice of the angular function $g$ produces a one-homogeneous function, $u^{g}$ say, defined by \eqref{ug}.  We shall employ the notation
\[u^{g}(x) = R g\left(\frac{x}{|x|}\right) \]
throughout the rest of the paper.

In the following, we denote the $m \times n$ real matrices by $\mathbb{R}^{m \times n}$, and unless stated otherwise we sum over repeated indices.  The elliptic system initially under consideration is
\begin{equation}
 \label{system} \frac{\partial}{\partial x_{q}}A_{pq}(x,\nabla u) = 0, \ \ \ 1 \leq p \leq m, \end{equation}
which is to be understood in the distributional sense, namely
\[\int_{B} A(\nabla u) \cdot \nabla \varphi \,dx = 0 \ \ \ \ \forall \varphi \in C_{c}^{1}(B, \mathbb{R}^{m}).\]
Here, $B$ is the unit ball in $\mathbb{R}^{2}$.

The function $A: B \times \mathbb{R}^{m \times 2} \to \mathbb{R}^{m \times 2}$ is $C^{1}$ and uniformly elliptic in the $F$ argument, that is
\begin{equation}
 \label{elliptic} \frac{ \partial A_{pq}(x,F)}{\partial F_{rs}}a_{p}b_{q}a_{r}b_{s} \geq \nu|a|^{2}|b|^{2} 
\end{equation}
for all $a \in \mathbb{R}^{m}$ and $b \in \mathbb{R}^{2}$.  The constant $\nu$ is independent of $x$ and $F$.  We say that a $C^{2}$ function $f: \mathbb{R}^{m \times n} \to \mathbb{R}$ is strongly rank-one convex if there is $\mu >0$ such that 
\[\frac{\partial^{2} f(F)}{\partial F_{rs} \partial F_{pq} } a_{r}b_{s}a_{p}b_{q} \geq \mu |a|^{2}|b|^{2}\]
for all $a \in \mathbb{R}^{m}$, $b \in \mathbb{R}^{n}$.  Replacing $a_{r}b_{s}$ and $a_{p}b_{q}$ respectively with $\pi_{rs}$ and $\pi_{pq}$ on the left, and $|a|^{2}|b|^{2}$ with $|\pi|^{2}$ on the right, leads to the usual definition  of strong convexity 
for $f: \mathbb{R}^{m \times n} \to \mathbb{R}$.

Other, standard notation includes $||\cdot ||_{k,p}$ for the norm on the Sobolev space $W^{k,p}$, $||\cdot||_{p}$ for the norm on $L^{p}$, and $\rightharpoonup$ to represent weak convergence in both of these spaces.    The tensor product of two vectors $a \in \mathbb{R}^{m}$ and $b \in \mathbb{R}^{n}$ is written $a \otimes b$; it is the $m \times n$ matrix whose $(i,j)$ entry is $a_{i}b_{j}$.   The inner product of two matrices $X,Y \in \mathbb{R}^{m \times n}$ is $X \cdot Y = \tr(X^{T}Y)$.   This obviously holds for vectors, too.   The $2 \times 2$ matrix $J$ will represent a rotation anticlockwise through $\frac{\pi}{2}$ radians, so that
\[J = \left(\begin{array}{c c } 0  & -1 \\ 1 & 0  \end{array} \right) .\]

\section{A regularity result for one-homogeneous maps}
We begin by establishing a regularity result for one-homogeneous stationary points of elliptic systems in two dimensions.  A differencing method in the angular variable is used to show that $W^{1,2}$ solutions must in fact be $W^{2,2}$.   It then follows from the Sobolev embedding theorem that all $W^{1,2}$ one-homogeneous stationary points must be Lipschitz, which observation is useful later in the paper.   Although we have focused on the case of $L^{2}$ integrable weak derivatives, it is plausible that similar arguments could be used to improve the regularity of $u^{g}$ in $W^{1,p}$ with $p \neq 2$, provided the growth and ellipticity hypotheses are suitably modified.  We do not do this since the improvement of regularity is ultimately put to a negative use in showing that such solutions are either affine or could not have existed in the first place.  

\begin{lemma}\label{firstlemma} Let $\om \subset \mathbb{R}$ be open and suppose it contains $0$.  Let $A: \om \times \mathbb{R}^{m \times 2} \to \mathbb{R}^{m \times 2}$ be $C^{1}$ in both arguments, and let $u^{g}(x)=|x|g(\theta)$ be in $W^{1, 2} (\om; \mathbb{R}^{m})$ and such that 
\begin{equation}\label{frog1} \int_{\om} A(x, \nabla u ^{g}(x)) \cdot \nabla \varphi(x) \, dx = 0 \end{equation}
for all $\varphi \in C^{1}_{c}(\om; \mathbb{R}^{m})$.  Suppose further that there are positive constants $C, \nu$ independent of $x$ such that 
\begin{itemize} \item[\textrm{(H1)}] $|D_{F}A(x,F)| \leq C(1+|F|)$ for all $F \in  \mathbb{R}^{m \times 2}$
\item[\textrm{(H2)}] $\frac{\partial A_{ij}}{\partial F_{rs}} (x,F) a_{i}b_{j}a_{r}b_{s} \geq \nu |a|^{2}|b|^{2}$ for all $a \in \mathbb{R}^{m}$ and $b \in \mathbb{R}^{2}$.   \end{itemize}
Then $g \in W^{2,2}(\mathbb{S}^{1}; \mathbb{R}^{m})$. 
\end{lemma}

\begin{remark}\label{matousek}\emph{The same result can be obtained if we replace the assumption that $0$ is an interior point of $\om$ with the assumption that $\om$ contains an (open) annulus.  The proof requires only minor changes to the last step of the argument given below.   }
\end{remark}
\begin{proof} Since $0 \in \om$ and $\om$ is open we assume without loss of generality that $\om = B(0,\delta)$ for sufficiently small $\delta$.  Let  
\[Q(h)= \left (\begin{array}{l l } \cos h & \sin h \\ -\sin h & \cos h  \end{array} \right)\]
for all real $h$, and let $B(x_{0}, 2\rho_{0}) \subset \om \setminus\{0\}$.   Let $\varphi$ be a smooth test function with support in $\om \setminus\{0\}$ and define
\[\varphi^{h}(x) = \varphi(Q(h)x)\]
for all real $h$.    Since $\varphi$ has compact support in $\om \setminus \{0\}$ then so does $\varphi^{h}$ for all sufficiently small $h$. Inserting $\varphi^{h}$ into \eqref{frog1} and changing variables we have 
\[\int_{\om}A(Q(h)x, \nabla u^{g} (Q(h)x))Q(h) \cdot \nabla \varphi(x) \,dx = 0.\] 
Recall that in polar coordinates one has, for non-zero $R$ and with $e_{R}=(\cos \theta, \sin \theta)$ and $e_{\theta}=(-\sin \theta, \cos \theta)$, that
\begin{equation}\label{altnab}\nabla \varphi =  \varphi_{,_{R}} \otimes e_{R} + \frac{1}{R}\varphi_{,_{\theta}}  \otimes e_{\theta}, \end{equation}
which gives, on setting 
\[A^{h}(x)= A(Q(h)x, \nabla u^{g} (Q(h)x)),\]
in the above that 
\[\int_{\om} A^{h}Q(h) \cdot \nabla \varphi\, dx = \int_{\om} (Ae_{R})^{h} \cdot \varphi_{,_{R}} + (Ae_{\theta})^{h}\cdot \frac{\varphi_{,_{\theta}}}{R} \,dx . \]
Here we have used the fact that $Q(h)e_{R}(\theta)=e_{R}(\theta + h)$, and similarly for $Q(h)e_{\theta}(\theta)$.   Hence, on using the notation 
\[\Delta_{h}z(x)=\frac{1}{h} (z(Q(h)x)-z(x))\]
for any function $z$ on $\mathbb{R}^{2}$ (which includes the matrix valued functions), 
\begin{equation}\label{difference1}
\int_{B(x_{0}, 2 \rho_{0})}\Delta_{h}(A(x,\nabla u^{g})e_{R}) \cdot \varphi_{,_{R}} + \frac{1}{R}\Delta_{h}(A(x,\nabla u^{g})e_{\theta}) \cdot \varphi_{,\theta}  \, dx = 0.
\end{equation}
Now,
\begin{eqnarray*} \Delta_{h}(A(x, \nabla u^{g})e_{R}) &= & A^{h}(x)\Delta_{h}(e_{R}) + \frac{1}{h}(A^{h}(x)-A(x,\nabla u^{g}(Q(h)x))e_{R} \\ 
                                                        & & + \frac{1}{h}(A(x,\nabla u^{g}(Q(h)x))-A(x,\nabla u^{g}(x)))e_{R}.
\end{eqnarray*}
Write this as 
\[\Delta_{h}(A(x, \nabla u^{g})e_{R}) = T_{1}(h)+T_{2}(h), \]
where 
\begin{eqnarray*}T_{1}(h) &= & A^{h}(x)\Delta_{h}(e_{R}) + \frac{1}{h}(A^{h}(x)-A(x,\nabla u^{g}(Q(h)x))e_{R} \\
T_{2}(h) & = & \frac{1}{h}(A(x,\nabla u^{g}(Q(h)x))-A(x,\nabla u^{g}(x)))e_{R}.\end{eqnarray*}
The differentiability and growth hypotheses on $A$ together with the assumptions on $g$ imply 
\[||T_{1}(h)||_{L^{2}(B(x_{0}, 2 \rho_{0}))} \leq c(||g||_{1,2}, x_{0}, \rho_{0}) \]
for some constant $c$ depending only on the quantities indicated (and in particular not on $h$).
A similar procedure can be followed for the other differenced term appearing in \eqref{difference1}.  The results are
\[\Delta_{h}(A(x, \nabla u^{g})e_{\theta}) = S_{1}(h)+S_{2}(h), \]
where
\[||S_{1}(h)||_{L^{2}(B(x_{0}, 2 \rho_{0}))} \leq c(||g||_{1,2}, x_{0}, \rho_{0}) \]
for some constant $c$ depending only on the quantities indicated (and in particular not on $h$), and
\[S_{2}(h)= \frac{1}{h}(A(x,\nabla u^{g}(Q(h)x))-A(x,\nabla u^{g}(x)))e_{\theta}.\]
As is usual in these cases, we write, for each  $1 \leq p \leq m$,
\[(S_{2}(h))_{p}= (\overline{DA})_{pqrs}(\Delta_{h}(\nabla u^{g}))_{rs}(e_{\theta})_{q},\]
where 
\[(\overline{DA})_{pqrs}= \int_{0}^{1}\frac{\partial A_{pq}}{\partial F_{rs}}(x, (1-t)\nabla u^{g}(x) + t \nabla u^{g} (Q(h)x)) \,dt.\]
Similarly, 
\[(T_{2}(h))_{p}= (\overline{DA})_{pqrs}(\Delta_{h}(\nabla u^{g}))_{rs}(e_{R})_{q}.\]

The rest of the proof consists in choosing $\varphi$ suitably and applying the ellipticity hypothesis to show that the quantity 
\[\int_{B(x_{0}, \rho_{0})}\frac{|\Delta_{h}g_{,_{\theta}}|^{2}}{R} \,dx \]
is bounded above independently of $h$.  One can then conclude the proof by applying Nirenberg's lemma.

Now
\[\Delta_{h}(\nabla u^{g})(x)= \Delta_{h}(g_{,_{\theta}}) \otimes e_{\theta} + g_{,_{\theta}}(\theta+h) \otimes \Delta_{h}(e_{\theta}) +\Delta_{h}(g \otimes e_{R}),\]
where the second and third terms are bounded in $L^{2}$ independently of $h$.  Let $\varphi= \eta^{2}\Delta_{h}g$, where $\eta$ is a smooth function with support in $B(x_{0}, 2\rho_{0})$ and satisfying $\eta = 1$ in $B(x_{0}, \rho_{0})$, $|\nabla \eta| \leq \frac{c}{\rho_{0}}$.  With this choice of $\varphi$ it can be checked that  
\begin{equation}\label{difference2}S_{1}(h) \cdot \frac{\varphi_{,_{\theta}}}{R} =   U_{1}(\eta, \nabla \eta, h, \Delta_{h}g, g_{,_{\theta}}) \cdot \frac{\eta \Delta_{h}(g_{,_{\theta}})}{R} + U_{2}(\eta, \nabla \eta, h, \Delta_{h}g, g_{,_{\theta}}),
\end{equation}
where $||U_{1}||_{2}$ and $||U_{2}||_{2}$ are bounded above independently of $h$.
Similarly, 
\begin{eqnarray*}S_{2}(h) \cdot \frac{\varphi_{,_{\theta}}}{R} & = &  \eta^{2}R (\overline{DA})_{pqrs}\frac{(\Delta_{h}g_{,_{\theta}})_{p}}{R}\frac{(\Delta_{h}g_{,_{\theta}})_{r}}{R}(e_{\theta})_{q}(e_{\theta})_{s} \\  & & + U_{3}(\eta, \nabla \eta, h, \Delta_{h}g) \cdot \frac{\eta \Delta_{h}(g_{,_{\theta}})}{R},
\end{eqnarray*}
where $||U_{3}||_{2}$ is bounded above independently of $h$.    The quantity $(T_{1}(h)+T_{2}(h))\cdot \varphi_{,_{R}} $ contributes only terms which appear on the right-hand side of \eqref{difference2}.  Therefore \eqref{difference1} can be written in  the form 
\begin{eqnarray}\label{difference3}
 \int_{B(x_{0}, 2 \rho_{0})}R (\overline{DA})_{pqrs}\frac{(\eta \Delta_{h}g_{_{\theta}})_{p}}{R}\frac{(\eta \Delta_{h} g_{,_{\theta}})_{r}}{R}(e_{\theta})_{q}(e_{\theta})_{s} \,dx & = & \nonumber \\  \ \ \ \ \ \ \ \int_{B(x_{0},2 \rho_{0})} \left( V_{1} \cdot \frac{\eta (\Delta_{h}g_{,_{\theta}})_{p}}{R}  +    V_{2} \right)\,dx,
\end{eqnarray} 
where the $V_{i}$ do not depend on $\Delta_{h}(g_{,_{\theta}})$ and $||V_{i}||_{2}$ are bounded above independently of $h$ for $i=1,2$.  
The uniform ellipticity of $A$ implies that 
\[\int_{B(x_{0}, 2 \rho_{0})}R (\overline{DA})_{pqrs}\frac{(\eta \Delta_{h}g_{,_{\theta}})_{p}}{R}\frac{(\eta \Delta_{h} g_{,_{\theta}})_{r}}{R}(e_{\theta})_{q}(e_{\theta})_{s} \,dx \geq  \nu \int_{B(x_{0},2 \rho_{0})} \frac{|\eta \Delta_{h} g_{,_{\theta}}|^{2}}{R} \,dx.\]
Using standard inequalities, the right-hand side of \eqref{difference3} can be bounded above by 
\[\frac{\eps^{2}}{2}\int_{B(x_{0}, 2 \rho_{0})} \frac{|\eta \Delta_{h} g_{,_{\theta}}|^{2}}{R} \,dx + \frac{1}{2 \eps^{2}} ||R^{-\frac{1}{2}}V_{1}||_{2}^{2}+||V_{2}||_{1},\]
where $\eps \neq 0$ may be chosen as small as we please, and in particular smaller than $\frac{\nu}{2}$.   Absorbing the term in $\frac{\eps^{2}}{2}$ into the term
\[\nu \int_{B(x_{0},2 \rho_{0})} \frac{|\eta \Delta_{h} (g_{,_{\theta}})|^{2}}{R} \,dx,\]
we have
\[\frac{\nu}{2} \int_{B(x_{0},2 \rho_{0})} \frac{|\eta \Delta_{h}(g_{,_{\theta}})|^{2}}{R} \,dx \leq \frac{1}{2 \eps^{2}} ||R^{-\frac{1}{2}}V_{1}||_{2}^{2}+||V_{2}||_{1} .\]
Recalling that $\eta=1$ on $B(x_{0}, \rho_{0})$, it follows that 
\begin{equation}\label{difference4}\int_{B(x_{0},\rho_{0})} |\Delta_{h}( g_{,_{\theta}})|^{2}\,dR \,d\theta \leq c, \end{equation}
where $c$ is independent of $h$.   Choosing $|x_{0}| = \frac{\delta}{2}$, $\rho_{0}=\frac{\delta}{8}$ and $\gamma$ to be the smallest postive solution of $\tan \gamma= 15^{-\frac{1}{2}}$, it can be checked that the subdomain
\[B_{\rho_{0}}:=\left\{x=R( \cos \theta, \sin \theta): \theta_{0}- \frac{\gamma}{2} \leq \theta \leq \theta_{0} + \frac{\gamma}{2}\right\} \cap B(x_{0}, \rho_{0}) \]
is such that 
\[\sup\{R: R(\cos \theta , \sin \theta) \in B_{\rho_{0}}\}- \inf\{R: R(\cos \theta , \sin \theta) \in B_{\rho_{0}}\} >  \rho_{0}\]
for each fixed $\theta \in [\theta_{0}-\frac{\gamma}{2}, \theta_{0}+ \frac{\gamma}{2}]$.  Therefore from \eqref{difference4},
\[\rho_{0} \int_{\theta_{0}-\gamma}^{\theta_{0}+\gamma} |\Delta_{h}( g_{,_{\theta}} )|^{2} \,d\theta \leq c.\]
Hence $g \in W^{2,2}_{\textrm{loc}}(\mathbb{S}^{1}, \mathbb{R}^{m})$.   Since $\gamma$ is independent of $\rho_{0}$ it follows that $g \in W^{2,2}(\mathbb{S}^{1}, \mathbb{R}^{m})$.
\end{proof}

The previous result can be used as follows:

\begin{theorem}\label{t1} Let $u$ be a $W^{1,2}(B, \mathbb{R}^{m})$ one-homogeneous solution to \eqref{system},   where $A$ satisfies
\begin{itemize}
 \item[\textrm{(H1)}]  
$A(x,F)$ is uniformly elliptic and $C^{1}$ in the gradient argument $F$; 
\item[\textrm{(H2)}] $|x|\partial_{x_{i}}A(x, F)$ is continuous on $(B\setminus\{0\}) \times \mathbb{R}^{n \times 2}$ for $i=1,2$;
\item[\textrm{(H3)}] $\lim_{R \to 0} R \partial_{x_{i}}A(x, \nabla u) = 0$ for $i=1,2$.\end{itemize}
Then $u$ is linear.
\end{theorem}

\begin{remark}\label{remark1}\emph{(H3) can hold for for functions $A$ whose spatial derivatives are singular at the origin, for example when 
\[c(F)R^{-\sigma} \leq \left|\frac{\partial A(x,F)}{\partial x_{i}}\right| \leq C(F)R^{-\sigma}\]
for $\sigma \in (0,1)$ and appropriate functions $c$ and $C$.}  \end{remark}

\begin{proof}Writing $u$ as 
\[u=Rg(\theta),\]
where 
\[g(\theta):=u(\cos \theta, \sin \theta),\]
it follows that $\nabla u$ depends only on the angular variable $\theta$:
\begin{equation}\label{gradu}\nabla u = g \otimes e_{R} + g' \otimes e_{\theta}.\end{equation}
As observed by Phillips, the elliptic system \eqref{system} can be written as
\begin{equation}\label{weakform}0 = \partial_{R} (A(x,\nabla u) e_{R})+ \frac{1}{R}\left(\partial_{\theta}(A(x,\nabla u)e_{\theta})+A(x,\nabla u) e_{R}\right). \end{equation}
By Lemma \ref{firstlemma}, we may assume that $u$ is Lipschitz.
From the independence of $\nabla u ^{g}$ on $R$ it follows that
\[R\partial_{R}(A(x,\nabla u^{g}) e_{R}) = R \frac{\partial A(x,\nabla u^{g})e_{R}}{\partial x_{i}} (e_{R})_{i}.\]
%m_{pq}(x)=R\frac{\partial A_{pq}}{\partial x_{i}}(x,\nabla u)(e_{\theta})_{i}(e_{\theta})_{q}\]
Since $g'$ is essentially bounded and by (H2) it follows that for each fixed $R >0$ the function  
\[\theta \mapsto R \partial_{R} (A(x,\nabla u) e_{R})+ A(x,\nabla u) e_{R}\]
is essentially bounded. Therefore from \eqref{weakform} for each fixed $R$ the function $\theta \mapsto \partial_{\theta}(A(x,\nabla u)e_{\theta})$ has a continuous representative.   Now we set about improving the regularity of the angular function $g(\theta)$ using the ellipticity hypothesis \eqref{elliptic}.   We refer the reader to \cite{Ph02} for a clear exposition in the $x-$independent case; in our case a similar argument works because \eqref{elliptic} is a uniform condition.

Now $A$ is strongly rank-one monotone, that is if $F-G=\xi \otimes \eta$ then 
\[(A_{pq}(F)-A_{pq}(G))(F-G)_{pq} \geq \nu |\xi|^{2}|\eta|^{2}.\]
By taking $F(\theta)=\nabla u(\theta)$ and $G(\theta, \varphi)=g(\theta)\otimes e_{R}(\theta)+g'(\varphi) \otimes e_{\theta}(\theta)$ for any fixed $\theta$ and $\varphi$ in $[0,2 \pi]$ it follows that
\[ F(\theta)-G(\theta, \varphi) = (g'(\theta)-g'(\varphi)) \otimes e_{\theta}(\theta),\]
and hence that
\begin{equation}\label{monot}|A(x(R,\theta),\nabla u(\theta))e_{\theta}(\theta)-A(x(R, \theta), G(\theta, \varphi)) e_{\theta}(\theta)| \geq \nu |g'(\theta) - g'(\varphi)|.\end{equation}
But $|A(x(R,\theta), G(\theta, \varphi))e_{\theta}(\theta)-A(x(R,\varphi), \nabla u(\varphi))e_{\theta}(\varphi)| \to 0$ as $\theta \to \varphi$, and since for fixed $R$ the function $A(x(R, \theta), \nabla u (\theta)e_{\theta}(\theta)$ is continuous in $\theta$ it must be that the left-hand side of \eqref{monot}
converges to zero as $\varphi \to \theta$.  Thus $g'$ is continuous.  But then \eqref{weakform} implies that $\partial_{\theta}(A(x,\nabla u)e_{\theta})$ is a $C^{1}$ function of $\theta$ for each fixed $R$.  It is now possible to follow Phillips's argument with only minor changes to deduce that $g''$ exists.  

From \eqref{weakform}, we see that for each $1 \leq p \leq m$ 
\begin{equation}\label{conclusion}R \frac{\partial A_{pq}(x,\nabla u)}{\partial x_{i}}\left((e_{R})_{q}(e_{R})_{i}+(e_{\theta})_{q}(e_{\theta})_{i}\right)+\frac{\partial A_{pq}(x,\nabla u)}{\partial F_{rs}}(g+g'')_{r}(e_{\theta})_{s}(e_{\theta})_{q}=0.\end{equation}

%Suppose that $g$ is not linear.  Then there is a set $\omega$ of positive measure in $[0,2 \pi]$  such that $g''(\theta)+g(\theta) \neq 0$ if $\theta \in \omega$.  

Multiplying \eqref{conclusion} by $(g''(\theta)+g(\theta))_{p}$ (and summing over $p$, therefore) one has, on applying the ellipticity hypothesis, that
\begin{equation}\label{greene}
 \nu|g(\theta)+g''(\theta)|^{2} \leq  
-R \frac{\partial A_{pq}(x,\nabla u)}{\partial x_{i}}\left((e_{R})_{q}(e_{R})_{i}+ (e_{\theta})_{q}(e_{\theta})_{i}\right)(g(\theta)+g''(\theta))_{p} .\end{equation}
Letting $R \to 0$ and applying (H3) forces the right-hand side to converge to zero.  Therefore $g$ is linear.
\end{proof}

%\begin{remark}\label{remark2} \emph{The same conclusion follows if we assume (H1), (H2$^{\prime}$) and (H3), where  
%\begin{itemize}\item[\textrm{(H2$^{\prime}$)}] 
%$|x|\partial_{x_{i}}A(x, F)$ is continuous on $\partial B(0,R) \times \mathbb{R}^{n \times 2}$ for $i=1,2$ and some $0 < R < 1$.  %\end{itemize}
%It is also clear from the proof that $g$ must be twice differentiable if we assume only (H1) and (H2$^{\prime}$): this will be of use in the next section.}
%\end{remark}

%The most obvious way to test this whether the hypotheses stated in tare necessary is to consider a class of elliptic systems of the form 
%\[A(x,F)= a(|x|)DW(F)\]
%where $a$ is a function chosen so that (H3) is violated and $W$ is strongly rank-one convex; this  last assumption ensures that $A$ is strongly elliptic in the sense of \eqref{elliptic} above.   We pursue this idea in the rest of the paper.

%We remark that searching specifically for a one-homogeneous solution to the elliptic system to searching for singular solutions in the general sense.  Indeed the work of Sv, Mu and Sz already shows that highly singular solutions to systems in divergence form exist; these solutions are not one-homogeneous.   
\section{Singular one-homogeneous stationary points}\label{francois}  In trying to prove the optimality of the results above it is natural to consider elliptic functionals which violate some of the conditions (H1) - (H3).    In this section we consider a specific functional which fails to be $C^{1}$ in its gradient variable, thereby violating two of the conditions of Theorem \ref{t1}.  The functional depends on the spatial variable as well as on the gradient of the competing functions.   Specifically, it is shown that there are non-trivial one-homogeneous solutions $u^{g}$ to a stationarity equation associated with the functional 
\[E(u) = \int_{\mathfrak{a}(R_{0},R_{1})} \frac{1}{|x|^{2}} W(\nabla u) \,dx.\]
Here, $0 < R_{0} < R_{1}$ and $\mathfrak{a}(R_{0}, R_{1})$ is the annulus in $\mathbb{R}^{2}$ centred on zero and with inner and outer radii $R_{0}$ and $R_{1}$ respectively.     The reason for this particular choice of integrand will be made clear later on.     The choice of an annular domain (as opposed to a ball) is forced on us.  To see why, note that
\[E(u^{g}) = \ln\left(\frac{R_{1}}{R_{0}}\right) \int_{\sone} W(\nabla u^{g}) \,d\theta.\]
Thus $R_{0} >0$ for all but the most trivial of problems where the integral over $\sone$ is zero.   But by cutting $0$ out of the domain we can no longer argue that the one-homogeneous functions are singular at $0$.   Something else has to be done to induce a singularity, which is the theme of the example discussed below.  

There are two strands to the argument of this section and the rest of the paper: one is concerned with existence, the other with ensuring that the solution, should it exist, is not linear.   Recall that Phillips's theorem \cite{Ph02} and Theorem \ref{t1} above conclude either that solutions are affine or do not exist at all.  We wish to avoid both possibilities.

We turn first to the question of existence.   
The stationarity condition referred to above is the so-called Equilibrium equation, meaning that 
\[\frac{d}{d\eps}\arrowvert_{\eps = 0 }E(u(x+\eps \Phi(x))) = 0 \ \ \ \forall \ \Phi \in C_{c}^{1}(\mathfrak{a}(R_{0}, R_{1}), \mathbb{R}^{2}).\]
In the $x-$independent case, and under suitable hypotheses, it is implied by the Euler-Lagrange equation.  See, for instance, \cite{BOP91}. The integrand $W$ mentioned above will be polyconvex and singular, after the fashion of the well-known stored-energy functions introduced by Ball in \cite{Ba77, Ba82}:    for $F \in \mathbb{R}^{2 \times 2}$, let  
\begin{equation}\label{w}W(F)=\frac{1}{2}|F|^{2} + h(\det F), \end{equation}
where  $h(t) =\infty$ if $t \leq 0$, and $h(t) \to \infty$ as $t \to 0+$.   The function $h$ is positive, $C^{2}$ and strongly convex on $(0,\infty)$.   One immediate consequence of this choice for $W$ is that the functional
\begin{equation}\label{ione} I(g) = \int_{\sone}  W(\nabla u^{g}) \,d\theta 
 \end{equation} 
is sequentially weakly lower semicontinuous with respect to weak convergence in $W^{1,2}(\sone, \mathbb{R}^{2})$.   This is a special case of the well-known results of Ball and Murat \cite{BM84}; we shall return to it in Proposition \ref{ball-murat} below.  

The second reason for choosing $W$ as above is to ensure that the minimizer $g$ of $I$ in some appropriate class $\scc$, say, is not the angular part of a linear map.  That is, we wish to prevent 
\begin{equation}\label{linearg} g(\theta) = T e_{R}(\theta),\end{equation}
where $T$ is a constant $2 \times 2$ matrix.   (This condition is necessary and sufficient for the corresponding mapping $u^{g}(x)=Tx$ to be linear.)
Clearly, the success or otherwise of this approach will also depend on the class of functions $\scc$ over which $I$ is minimized.   For now, suppose that any element $g$ of $\scc$ satisfies
\begin{itemize} \item[(i)]$I(g) < \infty$, and 
\item[(ii)] $g$ has a continuous representative which visits the origin in $\mathbb{R}^{2}$ at least once. \end{itemize}
If we suppose that some element $g$ of $\scc$ (not necessarily the minimizer) satisfies \eqref{linearg} then, by condition (ii), $\rank T \leq 1$.  But then $I(g) = \infty$ because $\det \nabla u^{g}$ is identically zero, contradicting (i).   Thus \emph{no} element $g$ of $\scc$ satisfying (i) and (ii) is such that $u^{g}$ is linear.   This condition differs to other, topological methods in ensuring that solutions are not linear.  For example, in the second half of their paper \cite{BOP91}, Bauman et al achieve the same goal essentially by restricting attention to a subclass of double-twist maps.   These methods do not seem to apply to problems involving one-homogeneous mappings.  

The Equilibrium equation arises in the context of singular integrands in \cite{Ba82} and \cite{BOP91}, and it is an appropriate starting point in the solution of problems in nonlinear elasticity theory.   In particular, it has been used as a stepping stone on the way to proving that certain of its solutions also solve the Euler-Lagrange equation.    

The one-homogeneous solution $u^{g}$ we construct will turn out to be singular in the sense that it is has an unbounded gradient on a half-line in $\mathbb{R}^{2}$.  It will be shown that this prevents $u^{g}$ from solving the Euler-Lagrange equation associated with $E$.     We investigate why this is so by giving a fairly detailed description of the manner in which the minimizer of $I$ in $\scc$ visits the origin in $\mathbb{R}^{2}$.   

%Finally, we discuss an example of a one-homogeneous solution to the full Euler-Lagrange equation associated with a functional similar to $E$.  The catch here is that the solution has degree two, which is in keeping with established methods for ensuring that solutions are not linear.   

%The weak limit of a sequence of one-hom functions is again one-hom, as can be seen by appealing to the formula
%\[\nabla u^{g} = g \otimes e_{R} + g_{,_{\theta}} \otimes e_{\theta} \]
%and noting that $u^{g^{(j)}}$ weakly converges in $W^{1,2}(B, \mathbb{R}^{2})$ if and only if $g^{(j)}$ weakly converges in $W^{1,2}(\mathbb{S}^{1}, \mathbb{R}^{2})$.

  % This is how the map $g$ alluded to above is obtained.  The drawback is that it does not necessarily satisfy the Euler-Lagrange equation as%sociated to the functional $I$.  Instead, it satisfies what is often referred to as the Equilibrium equations, which amounts to stationarity with respect to \emph{inner} variations of the candidate minimizer.    It turns out that $g$ is only Lipschitz in those connected components $\omega$ of $\mathbb{S}^{1}$ where 
%\[\textrm{ess} \inf_{\omega} \det \nabla u^{g} > 0.\]
%Such sets $\omega$ do not fill a.e. $\mathbb{S}^{1}$. 

\subsection{Necessary conditions satisfied by a minimizer of $I$}

In the following we shall assume $g: \sone \to \mathbb{R}^{2}$, with the convention that $g(\theta)=g(\cos \theta, \sin \theta)$, where
$\theta$ represents the polar angle.    Let $\varphi: \sone \to \mathbb{R}$ be a smooth function with compact support in $\sone$, and let us refer to such $\varphi$ as test functions.

Let $\C$ be a subset of functions in $W^{1,2}(\sone, \mathbb{R}^{2})$ with the properties
\begin{itemize}
\item[(P0)] $\scc$ is non-empty and closed with respect to weak convergence in $W^{1,2}(\sone, \mathbb{R}^{2})$;
\item[(P1)] $I(g) < \infty$ for all $g \in \C$; 
\item[(P2)] for each $g \in \C$ and each test function $\varphi$ there is $\eps_{0}>0$ such that $g^{\eps} \in \C$ for all $\eps \in (-\eps_{0}, \eps_{0})$,  where
\[g^{\eps}(\theta):=g(\theta+\eps \varphi(\theta));\]
\item[(P3)] for each $g \in \C$, each smooth $\varphi: \ \sone \to \mathbb{R}$ and each $\eps$, the condition 
\[\inf_{\theta \in \sone}\{|1+\eps \varphi(\theta)|\} \geq \frac{1}{2}\]
implies that $(1+\eps \varphi)g \in \C$.
\end{itemize}
  
 We remark that the test functions $\varphi$ need not have compact support in $\sone$ in order that integration by parts functions properly.  Instead, the periodicity of these functions suffices.
  
\begin{prop}\label{ball-murat} Let $I(g) = \int_{\sone} W(\nabla u^{g}) \,d\theta $, where $W$ is as per \eqref{w}, and let $\scc$ satisfy (P0) and (P1) above.  Then $\scc$ contains a global minimizer of $I$.\end{prop}
\begin{proof}  The direct method of the calculus of variations applies.  Any minimizing sequence $g_{j}$ is bounded uniformly in $W^{1,2}(\sone, \mathbb{R}^{2})$, which in view of the expression 
\[\nabla u^{g_{j}} = g_{j} \otimes e_{R} + {g_{j}}^{\prime} \otimes e_{\theta},\]
means that, for a subsequence, 
\[\nabla u^{g_{j}} \rightharpoonup \nabla u^{g}\]
for some $g \in \scc$.   By \cite[Theorem X]{BM84}, $I$ is lower semicontinuous with respect to sequential weak $W^{1,2}$ convergence.  Thus $g$ globally minimizes $I$.
 \end{proof}

Next, we derive two weak equations, \eqref{eqm} and \eqref{ii}, that the global minimizer of $I$ in $\scc$ must satisfy.  Notice that $\scc$ has so far only been described in terms of fairly generic properties.  In particular, we have not used any condition on the number of visits that curves in $\scc$ make to the origin, nor indeed any other kind of `boundary condition'.  % This observation will be useful in the very last section of the paper.

The calculations involved are non-trivial because of the singular integrand.  In proving \eqref{eqm} we follow the useful precedent in \cite{Ba82} and the subsequent paper \cite{BOP91}.

\begin{prop}\label{eqmequation}  Let $\C$ have properties (P0), (P1), (P2) and (P3) above.  Let $g$ be a global minimizer of $I$ in $\C$, where
\[I(g)=\int_{\sone}W(\nabla u^{g})\,d \theta\] 
and where $W$ is given in \eqref{w}.   In addition to the properties of $h$ assumed above, we suppose that there is a fixed and positive $s$ such that  $t^{s}h(t)$ and $t^{s + 1}h'(t)$ remain bounded as $t\to 0+$.   Then $g$ satisfies
\begin{eqnarray}\label{eqm}
\int_{\sone}\left(f(d) + \frac{1}{2}|g'|^{2}-\frac{1}{2}|g|^{2}\right)\varphi_{,_{\theta}}  \,d\theta &  = &  0 \\
\label{ii} \int_{\sone} \left( |g'|^{2}+|g|^{2}+2dh^{\prime}(d)\right)\varphi + (g' \cdot g ) \varphi_{,_{\theta}} \,d\theta & = & 0
\end{eqnarray} 
for all test functions $\varphi: \sone \to \mathbb{R}^{2}$.   Here, $f(t)=th'(t)-h(t)$ and $d=\det \nabla u^{g}$.   In addition, the quantity 
$f(d) + \frac{1}{2}|g'|^{2}-\frac{1}{2}|g|^{2} \in L^{1}(\sone)$.
 \end{prop}
\vspace{1mm}

\noindent \emph{Proof of \eqref{eqm}} For any $\varphi$ we may choose $\eps$ so small that 
\[z^{\eps}(\theta):=\theta+\eps \varphi(\theta)\]
is a diffeomorphism; we denote its inverse by $\psi^{\eps}$.   Taking $g^{\eps}(\theta) =  g(z^{\eps}(\theta))$  as above and applying (P2), we may suppose that $g^{\eps} \in \scc$ and hence that
   \begin{equation}\label{classicalstationary} \lim_{\eps \to  0} \frac{I(u^{g^{\eps}})- I(u^{g})}{\eps} = 0 \end{equation}
   whenever the limit on the left-hand side exists.   Let 
\begin{eqnarray*}d(\theta) & = & \det \nabla u^{g}(\theta) \\ d^{\eps}(\theta) &  = & \det \nabla u ^{g^{\eps}}(\theta). \end{eqnarray*}
Changing variables, we compute
\begin{eqnarray*} 
 \int_{\sone}\frac{1}{2}|\nabla u^{g^{\eps}}|^{2} \,d \theta & = & \int_{\sone} \frac{1}{2}\left(|g'(z)|^{2}(z'(\psi^{\eps}(z)))^{2}+|g(z)|^{2}\right)\, \frac{dz}{z'(\psi^{\eps}(z))} \\
\int_{\sone}h(d^{\eps}(\theta))\,d\theta & = &  \int_{\sone} h(z'(\psi^{\eps}(z)) d(z) ) \, \frac{dz}{z'(\psi^{\eps}(z))}. \end{eqnarray*}

We focus on the derivative of the second term with respect to $\epsilon$.  For brevity, let $x=z'(\psi^{\eps}(z))$.   Write
\[\int_{\sone}\frac{h(d^{\eps}(\theta))-h(d(\theta))}{\eps} \,d\theta = \int_{\sone}\frac{h(xd(z))-h(d(z))}{\eps x}\,dz + \int_{\sone}\frac{h(d(z))}{\eps}\left(\frac{1}{x}-1\right)\,dz.\]
 We take the limit as $\eps \to 0$ in the second integral by applying a version of the dominated convergence theorem and using the fact that 
 \begin{equation}\label{fact1}\lim_{\eps \to 0} \frac{1}{\eps}\left(\frac{1}{x}-1\right) = -\varphi^{\prime}(z).\end{equation} 
Here, $z$ is treated as the dummy variable: this is valid once $\eps$ has been set to zero, as can be seen by looking at the definition of the diffeomorphism $z^{\eps}$ given above.   The first integral can be rewritten, using Fubini's Theorem and other standard results, as
\begin{equation}\label{secondquotient} \int_{\sone}\frac{h(xd(z))-h(d(z))}{\eps x}\,dz  =  \int_{\sone} \int_{0}^{1}\frac{h^{\prime}(txd(z)+(1-t)d(z))d(z)(x-1)}{\eps x}\,dz \,dt.\end{equation}
Note that the assumptions on $h$ and $x$ mean that the pointwise bound
\[\frac{h^{\prime}(txd(z)+(1-t)d(z))d(z)(x-1)}{\eps x} \leq C h(d(z)) \]
holds for some positive constant $C$ depending on $\varphi$.  This is sufficient to validate 
\eqref{secondquotient} above.  It is now straightforward to take the limit as $\eps \to 0$ by applying the dominated convergence theorem, using \eqref{fact1} and noting that the pointwise limit of the integrand in \eqref{secondquotient} is independent of $t$.  The results are
\begin{eqnarray}\label{quotient1}\nonumber \lim_{\eps \to 0}
\int_{\sone}\frac{h(d^{\eps}(\theta))-h(d(\theta))}{\eps} \,d\theta & = & 
\int_{\sone}\left(d(\theta)h^{\prime}(d(\theta)) - h(d(\theta))\right)\varphi'(\theta) \, d\theta \\
& = & \int_{\sone} f(d(\theta)) \varphi'(\theta) \, d\theta,
\end{eqnarray} 
according to the definition of $f$ given in the statement of the proposition.

Arguing similarly, it can be shown that 
\begin{equation} \label{quotient2}\lim_{\eps \to 0}  \int_{\sone}\frac{|\nabla u^{g^{\eps}}|^{2}-|\nabla u^{g}|^{2}}{2 \eps} \, d\theta  = \int_{\sone}\left(\frac{|g'(\theta)|^{2}}{2}- 
\frac{|g(\theta)|^{2}}{2}\right)\varphi'(\theta) \, d\theta .
\end{equation}

Equation \eqref{eqm} now follows by combining equations \eqref{classicalstationary}, \eqref{quotient1} and \eqref{quotient2}. 
The last sentence of the statement is a straightforward consequence of the growth assumptions on $h$.

\vspace{2mm} 
\noindent \emph{Proof of \eqref{ii}} Let   $\varphi: \sone \to \mathbb{R}$ be an arbitrary smooth function and take $\eps$ so small that 
\begin{equation}\label{posdet}\inf_{\theta \in \sone}\{|1+\eps \varphi(\theta)|\} \geq \frac{1}{2}.\end{equation}
Let $g^{\eps,\varphi} = (1+\eps \varphi) g$.  Applying (P3), $g^{\eps, \varphi} \in \mathcal{C}$ for each fixed $\varphi$ and all sufficiently small $\eps$.   Since $g$ is a minimizer,
 \[\lim_{\eps \to 0} \frac{I(g^{\eps,\varphi}) - I(g)}{\eps} = 0\]
whenever the limit on the left-hand side exists.   But 
\begin{eqnarray*}\frac{I(g^{\eps,\varphi}) - I(g)}{\eps} & = &  \int_{\sone}\frac{\frac{1}{2}|(1+\eps \varphi)g'+\eps \varphi'g|^{2}-\frac{1}{2}|g'|^{2}}{\eps} \,d\theta \\
& & + \int_{\sone} \frac{\frac{1}{2}(1+\eps \varphi)^{2}|g|^{2}-\frac{1}{2}|g|^{2}}{\eps}\,d\theta \\ 
& & + \int_{\sone} \frac{h((1+\eps \varphi)^{2}d)- h(d)}{\eps} \,d\theta .
 \end{eqnarray*}
The limit as $\eps \to 0$ of the first two integrals on the right is easily seen to be
\begin{equation}\label{halfofit}\int_{\sone} \varphi (|g'|^{2}+|g|^{2})+(g' \cdot g) \varphi'\,d\theta. \end{equation}
The third integral can be dealt with by writing
\[\int_{\sone} \frac{h((1+\eps \varphi)^{2}d)- h(d)}{\eps} \,d\theta = \int_{\sone}\int_{0}^{1}h^{\prime}(d+t(2\eps \varphi +\eps^{2}\varphi^{2}))(2 \varphi d +\eps d \varphi^{2}) \,dt \,d\theta,\]
where we have implicitly used the estimate
\[|d h^{\prime}(d+t(2\eps \varphi +\eps^{2}\varphi^{2}))| \leq C h(d).\]
The constant $C$ depends on $\varphi$ but not on $\eps$ or $\theta$.   By Fubini's theorem and an appropriate version of the dominated convergence theorem we have
\[\lim_{\eps \to 0} \int_{\sone} \frac{h((1+\eps \varphi)^{2}d)- h(d)}{\eps} \,d\theta = \int_{\sone}2d h^{\prime}(d)\varphi \,d\theta, \]
which when added to \eqref{halfofit} gives \eqref{ii}.  This completes the proof. $\square$

\vspace{2mm}
It turns out that if $g$ satisfies the stationarity conditions \eqref{eqm} and \eqref{ii} then $u^{g}$ is automatically a solution of a fully two-dimensional equilibrium equation associated with the functional
\begin{equation*}E(u)= \int_{\mathfrak{a}(R_{0},R_{1})} \frac{W(\nabla u)}{|x|^{2}}\, dx. \end{equation*}
The equilibrium equation associated with $E$ is derived below, after which we show in Proposition \ref{twoimpliesone} that it is implied by \eqref{eqm} and \eqref{ii}.   For the sake of brevity we let $\mathfrak{a}=\mathfrak{a}(R_{0},R_{1})$ in the rest of the paper.

\begin{prop} Let $u$ be a $W^{1,2}(\mathfrak{a}; \mathbb{R}^{2})$ function such that $E(u) < \infty$.  Let $\Phi: \mathfrak{a}\to \mathbb{R}^{2}$ be an arbitrary function of compact support in $\mathfrak{a}$. 
Let $\eps_{0}$ be such that the function
\[z^{\eps}(x)=x+\eps \Phi(x)\]
is a diffeomorphism of $\mathfrak{a}$ whenever $|\eps| \leq \eps_{0}$.   Define the inner variation $u^{\eps}$ of $u$ by 
\[u^{\eps}(x) =  u(z^{\eps}(x)).\]
Suppose that
\[\lim_{\eps \to 0} \frac{E(u^{\eps})-E(u)}{\eps}  = 0.\]
Then $u$ satisfies 
\begin{equation}\label{fullgeneraleqn}\int_{\mathfrak{a}} 2W(\nabla u)\frac{x}{|x|^{4}} \cdot \Phi + \frac{1}{|x|^{2}}M(\nabla u) \cdot \nabla \Phi \,dx=0,\end{equation}
where
\[M(F)=F^{T}DW(F)-W(F)\1\]
for all $2 \times 2$ matrices $F$. (M is the so-called Energy-Momentum tensor.) When $u$ is a one-homogeneous map then the two-dimensional equilibrium equation \eqref{fullgeneraleqn} simplifies to
 \begin{equation}\label{full2DEqm}\int_{\mathfrak{a}} \left \{\left (M(\nabla u)e_{R}+2W(\nabla u)e_{R}\right) \cdot \Phi + M(\nabla u) e_{\theta} \cdot \Phi_{,_{\theta}}\right\} \, \frac{dR \,d\theta}{R^{2}} = 0.
\end{equation}
\end{prop}

\begin{proof}  It is customary to change variables in calculations involving inner variations.  To this end, for each $\eps$ in the range $(-\eps_{0}, \eps_{0})$ let the map $x^{\eps}$ be such that $x^{\eps} \circ z^{\eps}(x)=x$ for all $x \in \mathfrak{a}$.    Split $E(u^{\eps})-E(u)$ into three integrals as follows:
\begin{eqnarray*}\frac{E(u^{\eps})-E(u)}{\eps} & = &  \int_{\mathfrak{a}}\frac{1}{\eps}\left(\frac{1}{|x^{\eps}(z)|^{2}} - \frac{1}{|z|^{2}} \right ) W(\nabla u (z) \nabla z^{\eps} (x^{\eps}(z))\frac{dz}{\det \nabla z^{\eps}(x^{\eps}(z))} \\
&  +& \int_{\mathfrak{a}}\frac{1}{|z|^{2}}\left(\frac{W(\nabla u (z) \nabla z^{\eps} (x^{\eps}(z)))-W(\nabla u (z))}{\eps}\right)\frac{dz}{\det \nabla z^{\eps}(x^{\eps}(z))} \\
&  + &  \int_{\mathfrak{a}}\frac{1}{|z|^{2}}\frac{W(\nabla u (z))}{\eps} \left( \frac{1}{\det \nabla z^{\eps}(x^{\eps}(z))}-1\right) \,dz \\
& =: & I_{1}+I_{2}+I_{3}.
\end{eqnarray*}
It follows from the proof of \cite[Appendix]{BOP91a} that
\begin{equation}\label{firstlimit}\lim_{\eps \to 0}(I_{2}+I_{3}) = \int_{\mathfrak{a}}\frac{1}{|z|^{2}}M(\nabla u (z))\cdot \nabla \Phi (z) \,dz. \end{equation}

To calculate $\lim_{\eps \to 0} I_{1} $ it helps to note that 
\[x^{\eps} (z)= z - \eps \Phi(z) +o(\eps) \ \ \ \textrm{as} \ \eps \to 0, \]
which is a consequence of the choice made for the diffeomorphism $z^{\eps}$.  Therefore the limit
\[\lim_{\eps \to 0} \frac{1}{\eps}\left(\frac{1}{|x^{\eps}(z)|^{2}} - \frac{1}{|z|^{2}} \right)=  \frac{2z}{|z|^{4}} \cdot \Phi(z)\]
holds pointwise; applying this and a suitable convergence theorem to $I_{1}$ gives
\begin{equation}\label{prescient}\lim_{\eps \to 0} I_{1} = \int_{\mathfrak{a}} 2W(\nabla u) \frac{x}{|x|^{4}} \cdot \Phi \, dx.
\end{equation}
Adding \eqref{firstlimit} to \eqref{prescient} and setting the resultant expression equal to zero (using $ \lim_{\eps \to 0} (I_{1}+I_{2}+I_{3})=0$) gives the two-dimensional equilibrium equation \eqref{fullgeneraleqn}.

The calculation so far applies to any map, regardless of whether it is one-homogeneous or not.  
The final form \eqref{full2DEqm} of the equilibrium equation applies only to one-homogeneous maps; it can be reached by first noting that
\[\int_{\mathfrak{a}}  \frac{1}{|x|^{2}}M(\nabla u(x)) \cdot \nabla \Phi \,dx = \int_{\mathfrak{a}} \frac{1}{R} M(\nabla u) e_{R} \cdot \Phi_{,_{R}} + \frac{1}{R^{2}}M(\nabla u ) e_{\theta} \cdot \Phi_{,_{\theta}} \,dR \,d\theta,  \]
where we have used the expression \eqref{altnab} for $\nabla \Phi$ in polar coordinates.
Since  $M(\nabla u) e_{R}$ depends only on $\theta$, the first term can be integrated by parts with respect to $R$, thereby giving
\[\int_{\mathfrak{a}}\frac{1}{R} M(\nabla u) e_{R} \cdot \Phi_{,_{R}} \,dR \,d\theta = \int_{\mathfrak{a}} M(\nabla u) e_{R} \cdot \Phi \, \frac{dR \,d\theta}{R^{2}}.\]
Therefore
\[\int_{\mathfrak{a}} \frac{1}{|x|^{2}}M(\nabla u(x)) \cdot \nabla \Phi \,dx = \int_{\mathfrak{a}} \left\{M(\nabla u) e_{R} \cdot \Phi + M(\nabla u ) e_{\theta} \cdot \Phi_{,_{\theta}}\right\} \,\frac{dR \,d\theta}{R^{2}}.\]
Converting the remaining term in \eqref{fullgeneraleqn} into polar coordinates and combining with the above gives \eqref{full2DEqm}, as required.
\end{proof}

\begin{prop}\label{twoimpliesone} Let $g$ solve \eqref{eqm} and \eqref{ii} in a weak sense.  Then the one-homogeneous map $u^{g}$ satisfies 
\begin{equation}\label{2DEqm} \int_{\mathfrak{a}}\left\{ \left (M(\nabla u^{g})e_{R}+2W(\nabla u^{g})e_{R}\right) \cdot \Phi + M(\nabla u^{g}) e_{\theta} \cdot \Phi_{,_{\theta}}\right\} \, \frac{dR \,d\theta}{R^{2}} = 0.\end{equation}
Thus $u^{g}$ is a solution of the two-dimensional equilibrium equation associated with the functional $E(\cdot)$.
\end{prop}
\begin{proof}  Recalling that $M(F)=F^{T}DW(F) - W(F) \1$, we compute
\begin{eqnarray*} M(\nabla u^{g}) & = &  |g|^{2} e_{R} \otimes e_{R}  +  (g' \cdot g) \, (e_{R} \otimes e_{\theta} + e_{\theta} \otimes e_{R})  \\ & & + |g'|^{2} e_{\theta} \otimes e_{\theta} +  
  \left(f - \frac{1}{2}(|g|^{2} + |g'|^{2})\right)\1 .\end{eqnarray*}
Therefore
\begin{eqnarray} \label{A} M e_{R} & = &  \frac{1}{2}|g|^{2} e_{R} + (g' \cdot g) \, e_{\theta} + \left(f - \frac{1}{2}|g'|^{2})\right) e_{R} \\
\label{B} M e_{\theta} & = &  \frac{1}{2}|g'|^{2}e_{\theta} + (g' \cdot g) \, e_{R} + 
\left(f - \frac{1}{2}|g|^{2})\right)e_{\theta}, \end{eqnarray}
where $M=M(\nabla u^{g})$ for short.   Fix $R$ and suppress for now its appearance in $\varphi(\theta)$, where $\varphi(\theta) = \Phi(\theta,R) \cdot e_{R}$.  Then 
\begin{eqnarray*} \int_{\sone} Me_{\theta} \cdot \Phi_{,_{\theta}} \,d\theta & = &  \int_{\sone} \left(\frac{1}{2}|g'|^{2}+f-\frac{1}{2}|g|^{2}\right)\left((\Phi \cdot e_{\theta})_{\theta}+\varphi\right) \,d\theta \\ &  & + \int_{\sone} \left (g' \cdot g\right)\left(\Phi_{,_{\theta}} \cdot e_{R}\right)\,d\theta \\
& = & \int_{\sone} \left(\frac{1}{2}|g'|^{2}+f-\frac{1}{2}|g|^{2}\right)\varphi + (g' \cdot g)\left(\Phi_{,_{\theta}} \cdot e_{R} \right)\,d\theta, \end{eqnarray*} 
where we have used \eqref{eqm} to pass from one line to the next.
Since 
\begin{eqnarray*}\int_{\sone}\left (Me_{R}+2W(\nabla u^{g})e_{R}\right) \cdot \Phi \,d\theta &  = &   \int_{\sone}\left(\frac{3}{2}|g|^{2}+\frac{1}{2}|g'|^{2} + 2(h+f)\right) \varphi \,d\theta \\ & & + \int_{\sone}  (g' \cdot g)\Phi \cdot e_{\theta} \,d\theta,\end{eqnarray*}
it follows that 
{\small\begin{eqnarray*}
\int_{\sone}\left (Me_{R}+2W(\nabla u^{g})e_{R}\right) \cdot \Phi + Me_{\theta} \cdot \Phi_{,_{\theta}} \,d\theta & = & 
 \int_{\sone} \left(|g|^{2} + |g'|^{2}+2dh^{\prime}(d)\right)\varphi \,d\theta \\& +&   \int_{\sone}(g' \cdot g)\left(\Phi \cdot e_{\theta} +\Phi_{,_{\theta}} \cdot e_{R}\right) \,d\theta \\ 
& = & \int_{\sone}\left(|g|^{2} + |g'|^{2}+2dh^{\prime}(d)\right)\varphi \,d\theta\\ & & + \int_{\sone} (g' \cdot g) \varphi' \,d\theta \\
& = & 0
\end{eqnarray*}}
by \eqref{ii}.   Dividing both sides of this expression by $R^{2}$ and integrating over $R \in [R_{0},R_{1}]$ yields \eqref{2DEqm}. 
\end{proof}

\section{A class of curves which visit the origin in $\mathbb{R}^{2}$ at least once}\label{scarlatti}

In the following we identify each $g$ in $W^{1,2}(\sone, \mathbb{R}^{2})$ with its continuous representative (which exists by the Sobolev embedding theorem).  Let
\[ \scc = \{ g \in W^{1,2}(\sone, \mathbb{R}^{2}): \ I(g) < \infty, \  \exists \, \theta_{0} \in [0,2\pi] \ \textrm{s.t.} \ g(\theta_{0}) = 0\}. \]

\begin{figure}[ht]
\centering
\includegraphics[height=0.65\textwidth, width=0.65\textwidth]{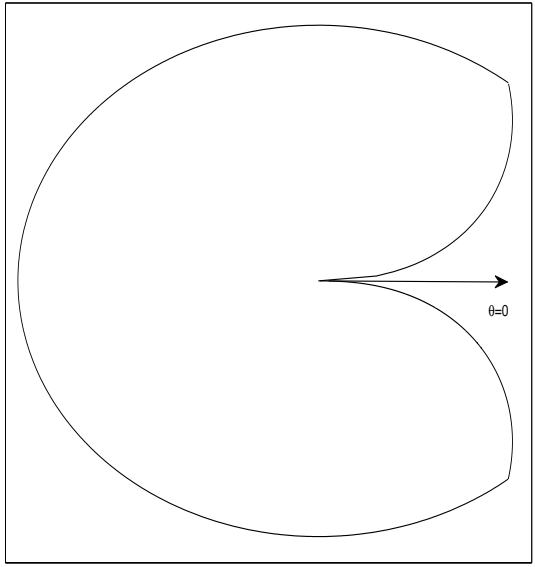}
\caption{An curve in the class $\scc$ corresponding to $s=\frac{1}{2}$, $l=\frac{1}{2}$ and $k=\frac{3}{8}$.}
%\psfrag{A}{$\theta = 0$}
%[position][psposition][scale][rotation]{LaTeX construction}
\label{cusp1}
\end{figure}

We prove that $\scc$ satisfies (P0)-(P3).   To verify (P0),  consider the map $g_{0}$ given in polar coordinates by 
\[g_{0}(\theta) = |\theta|^{k} e_{R}(\theta^{l}), \]
where $k$ and $l$ are constants to be chosen.   Let $\delta \in (0,\frac{1}{2})$ be fixed and let $\eta$ be a smooth, $2\pi-$periodic cut-off function satisfying $\eta(\theta) = 1$ in $|\theta | \leq \delta$ and with support in $[-2 \delta , 2 \delta]$.    
It is straightforward to check that the function
\[g:=\eta g_{0} + (1-\eta) e_{R}(\theta) \]
belongs to $\scc$ provided 
\[\int_{0}^{1} \theta^{2(k-1)} + \theta^{2(k+l)} + \theta^{-(2k+l)s} \,d\theta < \infty.\]
This holds if we choose $k = k_{\eps}$ and $l=l_{\eps}$, where
\begin{eqnarray*} k_{\eps}  & = &  \frac{1}{2}(1+ \eps)\\ 
l_{\eps} & = &  \frac{1}{s} - (1 + 2 \eps),\end{eqnarray*}
and $0 < \eps < \frac{2}{3s}$.   Thus $\scc$ is non-empty.  The argument of Proposition \eqref{ball-murat} shows that the weak limit $g$ of a sequence $\{g_{j}\} \subset  \scc$ has finite energy.  Moreover, since for each $g_{j}$ there is $\theta_{j}$ in $[0,2\pi]$ such that $g_{j}(\theta_{j})= 0 $, and since $g_{j}$ converges uniformly to $g$ on $[0,2\pi]$, it follows that $g$ visits the origin at least once.  Hence $\scc$ is weakly closed, and both (P0) and (P1) are verified.  Notice that if $g \in \scc$ then $\det \nabla u^{g}$ is strictly positive almost everywhere.  From the expression
\[\det \nabla u^{g} = Jg \cdot g',\]
and recalling that $J$ is the $2 \times 2$ matrix representing a rotation anticlockwise through $\frac{\pi}{2}$ radians, it can be inferred that the image of $g$ in $\mathbb{R}^{2}$ has a `handedness'.  For example, in Figure \ref{cusp1} below the curve is traversed anticlockwise.

To check (P2) note that a sufficiently small inner variation $g^{\eps}$, say, of any element $g \in \scc$ does not change the image set $g(\sone)$.  Thus any inner variation vists the origin at least once.  Furthemore, by a change of variables it can easily be seen that $I(g^{\eps}) < \infty$.   Hence $\scc$ is closed with respect to inner variations.   

To verify (P3), let  $\varphi: \sone \to \mathbb{R}$ be an arbitrary smooth function, take $\eps$ so small that 
\begin{equation}\inf_{\theta \in \sone}\{|1+\eps \varphi(\theta)|\} \geq \frac{1}{2}\end{equation}
and let $g^{\eps,\varphi} = (1+\eps \varphi) g$.
Then 
\begin{eqnarray*} \det \nabla u^{g^{\eps,\varphi}} & = &  Jg^{\eps,\varphi} \cdot (g^{\eps,\varphi})' \\
& = & (1+\eps \varphi)^{2} d(\theta) \\
& \geq & \frac{1}{4} d(\theta).
 \end{eqnarray*}
Using the growth conditions on $h$, it follows that $I(g^{\eps,\varphi}) < \infty$.   It is also clear that $g^{\eps,\varphi}$ visits the origin at least once because the same is true of $g$ by assumption.  Thus (P3) holds.

Having established (P0) - (P3), it follows that equations \eqref{eqm} and \eqref{ii} hold in the class $\scc$.  We begin their study with the observation that the Lipschitz solutions $g$ of \eqref{eqm} are those whose associated Jacobians $\det \nabla u^{g}$ are (positive and)  essentially bounded away from $0$. 

\begin{prop} \label{conformalgradient}  Let $g$ be a solution of \eqref{eqm}.  Let $\omega$ be a connected component of $\sone$, and let $\tilde \omega$ be a subset of $[0,2\pi]$ such that $\theta \in \tilde\omega$ if and only if $(\cos \theta, \sin \theta) \in \omega$.  
Then
 \[D_{\omega}:= \textrm{ess inf}\ \{d(\theta): \theta \in \tilde \omega \}\]
is strictly positive if and only if $g$ is Lipschitz on $\tilde \omega$.
%\begin{itemize}\item[(i)] $g$ is Lipschitz on $\tilde \omega$
%\item[(ii)] $g$ solves the weak Euler-Lagrange equation on $\tilde \omega$
%\item[(iii)] there is a conformal matrix $Q \in \RRtt$ such that 
%\[g(\theta)=Q e_{R}(\theta)\]
%for almost every $\theta$ in $\tilde \omega$.
%\end{itemize}
\end{prop}

\begin{proof}The weak form of \eqref{eqm} is 
\[\int_{\sone} \left(f(d)+\frac{|g'|^{2}}{2}-\frac{|g|^{2}}{2}\right)\varphi' \,d\theta = 0\]
 where, by Proposition \ref{eqmequation}, $f(d)+\frac{|g'|^{2}}{2}-\frac{|g|^{2}}{2} \in L^{1}(\sone)£$.  By DuBois-Reymond's lemma, there is $c \in \mathbb{R}$ such that 
\begin{equation}\label{DBReqm} f(d)+\frac{|g'|^{2}}{2}-\frac{|g|^{2}}{2} = c \ \textrm{a.e. in} \ \sone.\end{equation}
Recalling that $f(t)=th^{\prime}(t) - h(t)$, where $h$ is strictly convex, we see that $f^{\prime}(t)=th^{\prime \prime}(t) > 0$.  The hypotheses on $h$ ensure that $\lim_{t \to 0+} f(t) = -\infty$.  Therefore  
\[D_{\omega} > 0 \ \textrm{if and only if} \ \textrm{ess inf}\ \{f(d(\theta)): \theta \in \tilde{\omega}\}  > 0. \]
But then in view of \eqref{DBReqm} it follows that
\[\textrm{ess inf}\{f(d(\theta)): \theta \in \tilde{\omega}\}  > 0 \ \textrm{if and only if} \ \textrm{ess sup}\ \{|g'|: \theta \in \tilde{\omega}\} < \infty, \]
proving the claim.
\end{proof}

The stationarity conditions \eqref{eqm} and \eqref{ii} give some information about the regularity of an auxiliary quantity $z$ defined in \eqref{z} below.   We note that the function $z$ also appears naturally in \cite{BOP91} and \cite{YanPAMS}.

\begin{prop} \label{regularz} Let $g \in W^{1,2}(\sone, \mathbb{R}^{2})$ solve \eqref{eqm} and \eqref{ii}.  Then the function
\begin{equation}\label{z} z(\theta) : = \frac{1}{2}|g'(\theta)|^{2} + f ( d (\theta)) \end{equation}
lies in $W^{2,1}(\sone, \mathbb{R})$ and its weak derivatives satisfy
\begin{eqnarray}\label{firstderiv} z' & =  &  g \cdot g' \\
\label{secondderiv} z'' & =  &  2z + |g|^{2} + h(d).
\end{eqnarray}
In particular, $z$ is $C^{1}$. 
\end{prop}
\begin{proof} By \eqref{eqm},
 \[\int_{\sone}z \varphi' \,d\theta = \int_{\sone} -g \cdot g' \varphi \,d\theta.\]
Since $g \cdot g'$ is in $L^{1}$ it follows that the weak derivative of $z$ exists and satisfies \eqref{firstderiv}.  But \eqref{ii} implies that the weak derivative of $g \cdot g'$ exists and satisfies
\[ (g \cdot g')' = |g'|^{2} + 2dh^{\prime}(d) + |g|^{2}. \]
Rewriting this in terms of $z$ and $f$, and in view of \eqref{firstderiv}, it follows that $z''$ satisfies \eqref{secondderiv} as claimed. 
\end{proof}
 
\begin{remark}\emph{ This is as much as can be said while the regularity of $h \circ d$ is unknown.   We shall see later that there are circumstances in which the right-hand side of \eqref{secondderiv} becomes unbounded as a result of $d \to 0$.  In the same circumstances, Proposition \ref{regularz}  tells us that $z$ remains a $C^{1}$ map, even though $f(d) \to -\infty $ and $|g'| \to \infty$.}
\end{remark}

When \eqref{eqm} and \eqref{ii} hold we can say precisely when a solution is Lipschitz. 

\begin{prop}\label{singprop}  Let $g$ solve \eqref{eqm} and \eqref{ii}.  Then  
\[\textrm{ess inf}\, \{d(\theta): \theta \in \sone\} = 0\]  if and only if there is $\theta_{0}$ such that $g(\theta_{0})=0$.   In particular, by Proposition \ref{conformalgradient}, $g$ is a Lipschitz solution if and only if $g$ is never zero.
 \end{prop}
 
\begin{proof}  Assume that there is a sequence $\theta_{j} \to \theta_{0}$, say, such that $d(\theta_{j}) \to 0$.  It is claimed that $g(\theta_{0})=0$, where we identify $g$ with its continuous representative. Using \eqref{DBReqm} to identify
$\frac{1}{2}|g'|^{2} + f(d)$ with its continuous representative $c + \frac{1}{2}|g|^{2}$, and using the hypotheses on $h$ (and hence on $f$), it follows that
$|g'(\theta_{j})| \to \infty$.   Proposition \ref{regularz}  implies that $g \cdot g'$ is absolutely continuous, and hence in particular bounded (in modulus) on $\sone$ by some $C>0$.  Extracting a convergent subsequence from 
\[\psi_{j}:= \frac{g'(\theta_{j})}{|g'(\theta_{j})|},\] 
we may suppose $\psi_{j} \to \psi_{0}$, where $|\psi_{0}|=1$.  Therefore, since 
\[|g(\theta_{j}) \cdot \psi_{j}| \leq \frac{C}{|g'(\theta_{j})|}\]
it follows from the continuity of $g$ that
\[g(\theta_{0}) \cdot \psi_{0} = 0.\]
Also,  $d=Jg \cdot g'$ and $d(\theta_{j}) \to 0$ imply that
\[Jg(\theta_{0}) \cdot \psi_{0}=0.\]
Hence $g(\theta_{0})$ is orthogonal to both $\psi_{0}$ and $J \psi_{0}$, where $\psi_{0}$ is a unit vector, implying that $g(\theta_{0})=0$.

Suppose now that $g(\theta_{0})=0$ for some $\theta_{0}$.   If $\textrm{ess inf}\,\{d(\theta): \ \theta \in \sone\} > 0$ then we can argue as in Proposition \ref{conformalgradient} to conclude that $g$ is Lipschitz.  In particular, $\lim_{\theta \to \theta_{0}} Jg(\theta) \cdot g'(\theta)=0$, contradicting the assumption that $\textrm{ess inf}\,\{d(\theta): \ \theta \in \sone\} > 0$.  The last line in the statement of the Proposition is now straightforward.
\end{proof}

We now turn to the Euler-Lagrange equation associated to the functional $I$.     It is difficult to derive the equation by taking outer variations in the obvious way; indeed, we cannot be sure a priori that such a method is valid unless extra assumptions, such as those which feature in the lemma below, are made.  

\begin{lemma} \label{EL1d}  Let $g \in \scc$ be Lipschitz on the subinterval $\tilde{\omega}$ and suppose that $g$ minimizes $I$ in $\scc$.   Then for all $C^{1}$ functions $\xi : \sone \to \mathbb{R}^{2}$ with compact support in $\tilde{\omega}$, 
\begin{equation}\label{generalEL1D} \int_{\sone} DW(\nabla u^{g})e_{R} \cdot \xi + DW(\nabla u^{g})e_{\theta} \cdot \xi' \,d\theta = 0. \end{equation}
In particular, the equation  
\begin{equation}\label{EL1D}g -  h^{\prime}(d)Jg' - \left(g' + h^{\prime}(d) Jg \right)' = 0\end{equation}
holds in $\mathcal{D}^{\prime}(\tilde{\omega})$.    Moreover, $u^{g}$ satisfies a two-dimensional Euler-Lagrange equation
\begin{equation}\label{EL2D}\int_{\mathfrak{a}} \frac{1}{|x|^{2}}DW(\nabla u^{g}) \cdot \nabla \varphi \,dx = 0\end{equation}
among those $C^{1}$ functions $\varphi: \mathfrak{a} \to \mathbb{R}^{2}$ with the property that each map $\theta \mapsto \varphi(R,\theta)$ has support in $\tilde \omega$ for each fixed $R=|x|$ in $(R_{0},R_{1})$.   It follows by inspecting the proof of Theorem \ref{t1} that $g''$ exists in the classical sense on $\tilde \omega$ and is square integrable there.
\end{lemma}

\begin{proof}  Since $g$ is Lipschitz on $\tilde{\omega}$ it follows that both $d=g'(\theta) \cdot Jg(\theta)$ and $|g(\theta)|$ are bounded away from zero whenever $\theta \in \tilde{\omega}$.   Therefore for each $\xi$ as described in the statement of the lemma there is $\eps_{0} > 0 $ such that $|\eps| \leq \eps_{0}$ implies $g+\eps \xi \in \scc$.   A standard argument now implies \eqref{generalEL1D}.  
Inserting 
\begin{equation}\label{DW}DW(\nabla u^{g}) = g \otimes e_{R}+g' \otimes e_{\theta} +h^{\prime}(d)\left(Jg \otimes e_{\theta}-Jg' \otimes e_{R}\right) \end{equation}
into \eqref{generalEL1D} yields \eqref{EL1D}.  Equation \eqref{EL2D} exploits the one-homogeneity of $u^{g}$ in the same way as did the derivation of the two-dimensional equilibrium equation from \eqref{eqm} and \eqref{ii}.    The calculation in this case is straightforward and is left to the reader.
\end{proof}

\section{Spiral minimizers of $I$}\label{vivaldi}

We have fixed $\scc$ so that conditions (P0)-(P3) hold.  Therefore Propositions \ref{eqmequation},  \ref{twoimpliesone}, \ref{regularz} and \ref{singprop} all apply.   By definition, every element of $\scc$ has a zero.  Hence, in view of Proposition \ref{singprop}, none is Lipschitz.  It also follows from Proposition \ref{singprop} that the zeros of the global minimizer $g$ correspond exactly to the points where its gradient becomes unbounded.   Moreover, since $I(g) < \infty$, it follows that the set $\{\theta \in [0,2\pi]: \ d(\theta)=0\}$ is $\mathcal{H}^{1}-$null.   Therefore the set of zeros of $g$ is also $\mathcal{H}^{1}-$null.
This section is devoted to understanding the nature of the singularity in $g'$ associated with the zeros of $g$.     The following dichotomy is the starting point:
\begin{itemize}\item[(i)] either $g$ has no isolated zero, or  
\item[(ii)] there is $\theta_{0} \in [0,2\pi]$ and $\eps_{0} > 0$ such that $|g|>0$ on a relatively open interval in $[0,2\pi]$, one of whose endpoints is $\theta_{0}$.
\end{itemize}

If the zeros of $g$ were dense in $\sone$ then by the continuity of $g$, which, as before, can be inferred from Sobolev's embedding theorem, it would follow that $g$ is identically zero on $\sone$.   But $g=0$ is not a member of $\scc$, a contradiction.   Thus (i) is false, and at least one zero, $\theta_{0}$, say, of $g$ is isolated in the sense of (ii) above.  In consequence, $\nabla u^{g}$ has at least one line singularity which is mapped to $0$ under $u^{g}$.   It is tempting to conjecture that the minimizer $g$ has just one isolated zero, the reasoning being that it would be energetically unfavourable to incorporate more (thinking in terms of $I(g)$).   We do not pursue this conjecture here.    
%There is an obvious However, the argument in the proof of Proposition \ref{singprop} implies that the zeros of a minimizer correspond exactly to the points where $d$ vanishes.  Next, note that since the global minimizer of $I$ in $\scc$ has finite energy, it must be that       %Our intuition for restricting attention to isolated zeros is simply that multiple zeros are energetically costly, and as such are not likely to feature in an energy minimizer.   

Ultimately, we are interested in seeing whether $u^{g}$ solves the Euler-Lagrange equation associated with the functional
\[E(u)=\int_{\mathfrak{a}} \frac{1}{|x|^{2}}W(\nabla u) \,dx.\]
This can only be done once the behaviour of $g$ has been studied further.  Without loss of generality, we suppose that (ii) is true in a right neighbourhood of zero.   We also assume that $h(t)=t^{-s}$ for all positive $t$.  (The prescription $h(t)=+\infty$ for $t \leq 0$ continues to hold.)   This is not a restriction; it merely clarifies the subsequent analysis.  It will be shown that (ii) implies (ii'):
\begin{itemize}\item[(ii')] in a neighbourood of an isolated zero, $g$ rolls the annulus into an infinite spiral.\end{itemize}

\begin{prop}\label{bruce}Let $g$ be a minimizer of $I$ in $\scc$.   Assume that $|g(\theta)| > 0$ for $0 < \theta < \eps_{0}$ with $g(0)=0$.      Then, for $\theta \in (0,\eps_{0})$ and with 
\[h(t)=\left\{\begin{array}{ l l } t^{-s} & \textrm{if} \ t > 0 \\ +\infty & \textrm{otherwise,} \end{array} \right. \]
\begin{itemize} \item[(a)] a polar coordinate representation $g(\theta)=r(\theta)e_{R}(\gamma(\theta))$ is valid in $(0,\eps_{0})$, and in these coordinates the equations of stationarity \eqref{eqm} and \eqref{ii} become respectively
\begin{eqnarray}\label{1} \frac{1}{2}\left(r'^{2}+(rj)^{2}\right) & = & c + \frac{1}{2}r^{2} + (s+1)d^{-s} \\
\label{2} \frac{1}{4}(r^{2})^{\prime \prime} & = & c + r^{2} + d^{-s}, 
\end{eqnarray}
where $j=\gamma^{\prime}$;
\item[(b)] there is a nonnegative constant $\tau$ such that for $\theta \in (0,\eps_{0})$  
\begin{equation}\label{conservation1}sr^{2}(\theta)=\tau d^{s+1}(\theta)+ d^{s+2}(\theta), \end{equation}
where $d(\theta)=\det \nabla u^{g}$; 
\item[(c)] d is monotone increasing in a (right) neighbourhood of $\theta=0$.
\end{itemize}
\end{prop}

\begin{proof} The assumptions on $g$ are such that Lemma \ref{EL1d} applies with $\tilde{\omega}=(0,\eps_{0})$.    In particular, \eqref{EL1D} holds with $g''$ classically second differentiable at all points in $(0,\eps_{0})$.   This improvement in regularity means that the polar coordinate representation of $g$ on the interval $(0,\eps)$ is indeed valid.   Equations \eqref{1} and \eqref{2}
now follow from this and the concrete choice for $h$ made above.  

To prove (b), we take the inner product of  \eqref{EL1D} with $Jg$.  (We remark that \eqref{1} and \eqref{2} can be recovered from \eqref{EL1D} by taking its inner product with $g$ and $g'$ respectively.)   It helps to recall that $d=g^{\prime} \cdot Jg$, where $J$ is the rotation through $\pi/2$ radians anticlockwise, in order to see that $d'=g'' \cdot Jg$.   Also, $Jg \cdot Jg' = g \cdot g'= rr'$.  Therefore on all compact subintervals of $(0,\eps_{0})$
\[2rr' h'(d) + d' + h''(d)r^{2} d' = 0 .\]
It follows by integration that there is $-\tau$ such that 
\[d+r^{2}h'(d)=-\tau.\]
Rearranging this and inserting $h$ as described above yields
\[sr^{2}=d^{s+2}+\tau d^{s+1}\]
on $(0,\eps_{0})$, which is \eqref{conservation1}.   If $\tau \neq 0$ then its sign can be deduced as follows.   The improved regularity of $g$ demonstrated in part (a), together with the observation that $g=0$ if and only if $d = 0$, implies that $d$ is continuous on $(0,\eps_{0})$ and satisfies $d(\theta) \to 0$ as $\theta \to 0+$.  Therefore $\tau d^{s+1}$ dominates the right-hand side of \eqref{conservation1} as $\theta \to 0$, and it follows easily that $\tau > 0$ if it is non-zero.  

To prove (c) it suffices to show that ${r^{2}}^{\prime}> 0 $ near zero; one then appeals to \eqref{conservation1} to conclude that $d$ must also be strictly increasing.  Equation \eqref{2} implies that for $\theta$ sufficiently small and positive we may assume that $r^{2}$ is strongly convex.   Now $r^{2}(0)=0$ by hypothesis.  Translating the final line in the statement of Proposition \ref{regularz} into polar coordinates, we see that ${r^{2}}^{\prime}$ is continuous on all of $\sone$.  In particular,  if ${r^{2}}^{\prime}(0)$ were non-zero then it would imply that $r^{2}(\theta) < 0$ in either a left or right neighbourhood of $\theta= 0$; either way this is a contradiction.  Therefore ${r^{2}}^{\prime}(0)=0$, and hence by the strong convexity of $r^{2}$ it must be that ${r^{2}}^{\prime} > 0 $ on $(0,\eps)$ for some $\eps > 0$.   This concludes the proof.
\end{proof}

 It can be checked that when $\tau$ is strictly positive the solution curve winds only finitely many times around the origin, with smaller values of $\tau$ corresponding to higher winding numbers.     It therefore seems quite natural that the solution in the case $\tau=0$ is an infinite spiral.   However, we are not free to choose $\tau$: its value is imposed on us by the minimization process.  I cannot rule out the possibility that there are stationary points whose corresponding value of $\tau$ is strictly positive, but their existence is not proven by the methods used in this paper.   Instead, we focus below on showing that $\tau$ must be zero when it satisfies \eqref{conservation1} and when $g$ minimizes $I$ in $\scc$.

\begin{prop}\label{untimely}Let $r$ and $j$ be as in Proposition \ref{bruce} above, and recall in particular that they correspond to a minimizer $g$ of $I$ in $\scc$.    Suppose that the nonnegative constant $\tau$ satisfies  
\begin{equation}\label{tau} sr^{2}=\tau d^{s+1} + d^{s+2} \end{equation}
on $(0,\eps_{0})$.    Then $\tau=0$. \end{prop}

\begin{proof}  Assume for a contradiction that $\tau > 0$.   
It will be shown that there are variations of $\hat{g}$ in $\scc$ which lower the energy $I$, from which the result follows immediately.  
The proof is divided into $4$ parts, the first of which establishes some basic facts about $r$ and $j$.\vspace{2mm}

\textbf{Step 1.}
Let $\tilde{\tau}= \tau + d$ and note that \eqref{conservation1} implies 
\begin{equation}\label{r}r = \left(\frac{s}{\tilde{\tau}}\right)^{\frac{1}{2s}}j^{-\left(\frac{s+1}{2s}\right)}.\end{equation}
($r$ is still only defined implicitly by this expression.)  Note that  $\ttau > 0$ for all $\theta$ because $d >0$ a.e. and $\tau > 0$ by assumption. 
Since $d= r^{2}j$, it follows that 
\begin{equation}\label{h}d^{-s} = \frac{\tilde{\tau}j}{s}, \end{equation}
and hence from \eqref{1} that 
\begin{equation}\label{rprime} r' = \left(2\left(1+\frac{1}{s}\right) \tilde{\tau} j \right)^{\frac{1}{2}} F,\end{equation}
where 
\[F^{2} = 1 + \frac{r^{2}+2c}{2(1+\frac{1}{s}) \ttau j} - \frac{s^{\frac{1}{s}}}{2(1+\frac{1}{s})} \tilde{\tau}^{-(1+\frac{1}{s})}j^{-\frac{1}{s}}. \]
Differentiating \eqref{conservation1} with respect to $\theta$ gives 
\[2srr^{\prime} = ((s+1)\tau + (s+2) d) d^{s} d^{\prime}, \]
which on using \eqref{r}, \eqref{h} and \eqref{rprime} above gives
\begin{equation}\label{schubert} d^{\prime} = \frac{2 s^{\frac{1}{2s}} \ttau^{\frac{3}{2}-\frac{1}{2s}} \left(2\left(1+\frac{1}{s}\right) \right)^{\frac{1}{2}} j^{1-\frac{1}{2s}}F}{(s+1)\tau + (s+2) d}.\end{equation} 
But $d^{\prime} = 2 r r^{\prime}j + r^{2} j^{\prime}$, which, on eliminating $r, r^{\prime}$ and $d^{\prime}$ using the expressions given so far, shows that $j$ satisfies the equation
\begin{equation}\label{Yfirst}j^{\prime} + Y j^{2+\frac{1}{2s}} = 0 \end{equation}
on some interval $(0,\eps)$, where
\begin{equation}\label{Y}Y = 2 s^{-\frac{1}{2s}}\left(2\left(1+\frac{1}{s}\right) \right)^{\frac{1}{2}} F \ttau^{\frac{1}{2}+\frac{1}{2s}} \left(1 - \frac{\ttau}{((s+1)\tau + (s+2) d)}\right). \end{equation}
Since $d \to 0$ monotonically (by Propositon \ref{bruce}, part (c)), and since \eqref{h} holds, it follows that $j \to \infty$ monotonically as $\theta \to 0$.   This boundary condition allows us to solve, at least in principle, the differential equation \eqref{Y}.
It also follows from this that $F$ is very close to $1$ for all sufficiently small $\theta$.   This observation will be used below.\vspace{2mm}

\textbf{Step 2.}  Recall that the polar coordinate representation of the minimizing map is only known to be valid on some inerval $(0,\eps)$. Therefore its energy is represented by 
\[I(g) = \int_{0}^{\eps} \frac{1}{2}(r^{2}+{r^{\prime}}^{2}+(rj)^{2}) + d^{-s} \ d\theta + \int_{\eps}^{2\pi}W(\nabla u ^{g}) \,d\theta.\]
The first integral on the right-hand side can be written in terms of $j$ and $\ttau$ using \eqref{r},\eqref{h} and
 \eqref{rprime}.  The result is
 \[I(g) = \int_{0}^{\eps} ((s+1)F^{2}+1)\frac{\ttau j}{s}+ \frac{1}{2}s^{\frac{1}{s}}\ttau^{-\frac{1}{s}}j^{-(\frac{s+1}{s})}(1+j^{2}) \ d\theta + \int_{\eps}^{2\pi}W(\nabla u ^{g}) \,d\theta.\]
Since $\tau$ is fixed and $j \to \infty$ as $\theta \to 0$, the integrand of the first integral on the right is dominated by the term in $\ttau j$ as $\theta \to 0$.    One can infer from this, albeit informally, that a slightly smaller value of $\tau$ would suffice to lower the energy.  In practice, one has to be careful about changing $\tau$: the effect might be global, possibly even resulting in an overall increase in the energy.  In Step 3 below we vary $\ttau$ near zero whilst retaining its limiting value of $\tau$, thereby keeping the effect of the change local and hence controllable.  \vspace{2mm}

\textbf{Step 3.}  Let us define a variation $\hat{g}$ about $g$ in terms of the angular velocity $j$ and radial component $r$ of $g$ as follows.    Firstly, let $T= \eta \ttau$, where $\eta$ is a smooth map with support in $[\delta_{1}, \delta_{2}] \subset (0,\eps)$.  For now we think of $\eta$ as being close to $1$ in value; in this sense $\hat{g}$ is considered a perturbation of $g$.  The parameters $\delta_{1}$ and $\delta_{2}$ will be chosen shortly.   Define $\hat{r}$ by 
\[s{\hat{r}}^{2} = T d^{s+1},\]
where $d$ is as per \eqref{h} above and $\ttau = \tau + d$ as before.   Using the relation $\hat{d} = \hat{r}^{2}j$, it follows that 
\[\hat{d}^{-s} = \eta^{-s} d^{-s}. \]   
Therefore 
\begin{equation}\label{diffds}d^{-s} - \hat{d}^{-s} = (1- \eta^{-s})\frac{\ttau j}{s},\end{equation}
where we have used the expression \eqref{h} for $d$ from the previous step.  

Also, using \eqref{r} in conjunction with the definition of $\hat{r}$ given above it can be seen that
\begin{equation} \label{faure} \frac{1}{2}\left(r^{2} + (rj)^{2} - \hat{r}^{2} - (\hat{r}j)^{2}\right) = \frac{1}{2}(1-\eta) s^{\frac{1}{s}} \ttau^{-\frac{1}{s}} j^{-\frac{s+1}{s}}(1+j^{2}). \end{equation}   The only term it remains to compute in the integrand of
\[\int_{0}^{\eps} W(\nabla u^{g}) - W(\nabla u^{\hat{g}}) \,d\theta \]
involves $\frac{1}{2}({r^{\prime}}^{2} - (\hat{r}')^{2})$.    In what follows it will be convenient to combine this difference with the difference
$\hat{d}^{-s}- d^{-s}$.   Therefore we let
\[\Xi = \frac{1}{2}({r^{\prime}}^{2} - (\hat{r}')^{2}) + d^{-s}-\hat{d}^{-s}.\]
The next and final step of the proof analyses the behaviour of $\Xi$ for small values of $\theta$.   

\vspace{2mm}
\textbf{Step 4.}   Rewrite \eqref{schubert} as $d' = D j ^{1-\frac{1}{2s}}$, where
\begin{equation}\label{bigD} D = \frac{2 s^{\frac{1}{2s}} \ttau^{\frac{3}{2}-\frac{1}{2s}} \left(2\left(1+\frac{1}{s}\right) \right)^{\frac{1}{2}}}{(s+1)\tau + (s+2) d}.\end{equation}
Using this shorthand when differentiating 
\[\hat{r} = \left(\frac{T}{s}\right)^{\frac{1}{2}}d^{\frac{s+1}{2}},\]
using the definition of $T$ and the fact that $\ttau'=d'$, it can be seen that
\begin{eqnarray*}(\hat{r}')^{2} & =  & \frac{(s+1)^{2}}{4s}\ttau \eta d^{s-1}D^{2}j^{2-\frac{1}{s}} + \frac{(s+1)}{2}\eta D^{2} \ttau^{-1}j^{1-\frac{1}{s}} + \\ &  + &  \frac{s^{\frac{1}{s}}}{4}\eta^{2} D^{2} \ttau^{-(1+\frac{1}{s})}j^{1-\frac{2}{s}} + O(\eta').\end{eqnarray*} 
Recalling that $\eta$ remains close to $1$, we can further suppose that $\eta'$ is small.  (We will later choose $\eta$ so that this is the case.) Next, we form $\Xi$ by grouping together terms in a suitable way.  To make it explicit we first set
\[\Xi= A_{1} + A_{2}\]
where
\begin{eqnarray*} A_{1}:= \frac{1}{2}r'^{2} +d^{-s} - \hat{d}^{-s} - \frac{(s+1)^{2}}{8s}\ttau \eta d^{s-1}D^{2}j^{2-\frac{1}{s}} \\
 A_{2}: = - \frac{(s+1)}{4}\eta D^{2} \ttau^{-1}j^{1-\frac{1}{s}} -  \frac{s^{\frac{1}{s}}}{8}\eta^{2} D^{2}\ttau^{-(1+\frac{1}{s})}j^{1-\frac{2}{s}} + O(\eta').
\end{eqnarray*}
Replacing $D^{2}$ in $A_{1}$ using \eqref{bigD}, and using equations \eqref{rprime} and \eqref{diffds}, it follows that
\[A_{1} = \frac{\ttau j }{s}((1+s)F^{2}- s)(1-\eta)+\frac{(1+s)\ttau j}{s}F^{2} \eta \left(1- \left( \frac{\tau +d}{\tau + \frac{(s+2)}{(s+1)}d}\right)^{2} \right)  + o(1-\eta).\]
Combining this expression with the first two terms of $A_{2}$ and simplifying yields
\[\Xi = \frac{\ttau j }{s}((1+s)F^{2}- s)(1-\eta) + \frac{\ttau \eta j d^{2} F^{2}}{s}\frac{1-\eta \ttau}{\left(\tau + \left(\frac{s+2}{s+1}\right)d\right)^{2}} + o(1-\eta) + O(\eta').\]
Since $F^{2}$ converges monotonically to $1$ there exists $\theta_{0} > 0$ such that 
\[(1+s)F^{2}- s > \frac{1}{2} \ \textrm{if} \ \theta \in (0,\theta_{0}).\]
In particular, if $\eta < 1$ is enforced then the first two terms of $\Xi$ remain positive on $(0,\theta_{0})$.
Let $\psi$ be a fixed smooth function with support compactly contained in $(0,\theta_{0})$ and which satisfies $0 \leq \psi \leq 1$.  Let $\sigma$ be a small, positive parameter and let $\eta = 1 - \sigma \psi$.  Insert this choice of $\eta$ into the last expression for $\Xi$ and call the result $\Xi_{\sigma}$.  Then, by integrating $\Xi_{\sigma}$  over $(0,\eps)$ and applying a version of the dominated convergence theorem, we see that 
\[\int_{0}^{\eps} \Xi_{\sigma} \,d\theta > 0\] 
provided $\sigma$ is sufficiently small.   Finally, we have to include the terms 
\[\frac{1}{2}\left(r^{2}-\hat{r}^{2}+(rj)^{2} - (\hat{r}j)^{2}\right).\] 
But by \eqref{faure}, this quantity is strictly positive whenever $\eta < 1$, and is zero otherwise.   We conclude then that  $I(g)-I(\hat{g}) > 0 $, provided the perturbation $\hat{g}$ is defined as per Step 3.  This contradicts our assumption that $g$ minimizes $I$ in $\scc$.  

\end{proof}

\begin{prop}\label{hubbard}  Let $r$ and $j$ be as in Proposition \ref{bruce} above and let $\tau = 0$, so that 
\begin{equation}\label{conservation2} sr^{2}=d^{s+2} \ \textrm{on} \ (0,\eps_{0}). \end{equation}
Define 
\[ n(s) =   \frac{s+2}{2(s+1)} \]
and
\[ V = \left(1 + \frac{ d^{s}(r^{2}+2c)}{(s+2)}\right)^{\frac{1}{2}}. \]
Then $j$ satisfies
\begin{equation}\label{Ysecond}\frac{2(s+1)}{(s(2+s))^{\frac{1}{2}}}Vj^{2}+ j^{\prime} = 0 \ \textrm{on} \ (0,\eps_{0}), \end{equation}
subject to $\lim_{\theta \to 0+} j = \infty$, 
and 
\[r = s^{\frac{1}{2(s+1)}}j^{-n(s)}. \]
\end{prop}
\begin{proof}  The expression for $r$ given above follows directly from \eqref{conservation2}.   The differential equation \eqref{Yfirst}, which was found during the proof of the previous proposition, yields \eqref{Ysecond}.  The various coefficients and exponents can be evaluated using $\tau = 0$, $\ttau=d$ and $d=s^{\frac{1}{s+1}}j^{-\frac{1}{s+1}}$.  Integrating \eqref{Ysecond}, with 
\[\beta : = \frac{2(s+1)}{(s(2+s))^{\frac{1}{2}}},\]
yields, for $0 < \theta < \eps$, 
\[j^{-1}(\theta) = \int_{0}^{\theta} \beta V(\bar{\theta}) \,d\bar{\theta}.\] 
It follows from this and the expression given for $V$ above that there are positive constants $m < M$ such that
\begin{equation}\label{jbehaviour}\frac{1}{M\theta} \leq j(\theta) \leq \frac{1}{m\theta} \end{equation}
for $0 < \theta < \eps$.
Recalling that $j = \gamma'$, where $g=r(\theta)e_{R}(\gamma(\theta))$, and integrating this 
expression over $(0,\eps)$ yields 
\[\frac{1}{M} \ln \left(\frac{\theta}{\eps}\right) + \gamma(\eps) \leq \gamma(\theta) \leq   \frac{1}{m} \ln \left(\frac{\theta}{\eps}\right) + \gamma(\eps),\]
which is valid on $(0,\eps)$.   It also follows from the expression given above for $r$ and from \eqref{jbehaviour} that $r$ is bounded above and below by an expression of the form $C\theta^{n(s)}$, where $C$ is constant.  
\end{proof}

\subsection{Non-satisfaction of the Euler-Lagrange equation near an isolated singularity}

The Euler-Lagrange equation associated with $E$ is 
\[\int_{\mathfrak{a}} \frac{1}{|x|^{2}}DW(\nabla u^{g}) \cdot \nabla \Phi \,dx = 0 \ \ \forall \  \Phi \in C^{1}_{c}(\mathfrak{a}; \mathbb{R}^{2}). \]
In polar coordinates, this is
\begin{equation}\label{couperin}\int_{\mathfrak{a}} \frac{1}{R^{2}}(DW(\nabla u^{g})e_{R} \cdot \Phi_{,_{R}} + \frac{1}{R} DW(\nabla u^{g})e_{\theta} \cdot \Phi_{,_{\theta}})\,dx = 0.\end{equation}
From \eqref{DW}, it follows in particular that
\[DW(\nabla u^{g})e_{R} = g + sd^{-(s+1)} Jg'. \]
Now, $sd^{-(s+1)} = j$ and $Jg' = r' e_{\theta} -rj e_{R}$; hence,
\[DW(\nabla u^{g})e_{R} = r(1-j^{2})e_{R} + r'j e_{\theta}.\]
The asymptotic behaviours of $r$ and $j$ as $\theta \to 0$ are given in Proposition \eqref{hubbard} above; they imply that both the terms
$rj^{2}$ and $r'j$ are of order $\theta^{(n(s)-2)}= \theta^{-1-\frac{s}{2(s+1)}}$ as $\theta \to 0$, which is not $L^{1}$-integrable.  It is now easy to choose a test function $\Phi$ such that the Euler-Lagrange equation \eqref{couperin} fails.

\vspace{3mm}
\noindent \textbf{Acknowledgement}  I would like to thank Prof. Kewei Zhang and Dr Elaine Crooks for their helpful comments on this work, which was completed with the support of an RCUK Academic Fellowship.

\end{document}